\input amstex
\documentstyle{amsppt}

\tolerance 3000 \pagewidth{5.5in} \vsize7.0in
\magnification=\magstep1 \widestnumber \key{AAAAAAAAAAAAAAAAAA}
\NoRunningHeads \loadmsam \topmatter
\title On combinatorial
complexity of convex sequences
\endtitle
\author A. Iosevich, M. Rudnev and V. Ten
\endauthor
\address Alex Iosevich. Department of Mathematics,
The University of Missouri, Columbia, MO 652111, U.S.A.
\endaddress

\email iosevich\@math.missouri.edu \newline
http://www.math.missouri.edu/$\thicksim$iosevich/\endemail

\address Mischa Rudnev. School of Mathematics,
University Walk, Bristol BS8 1TW, U.K.
\endaddress
\email m.rudnev\@bris.ac.uk \endemail
\address Vladimir Ten. School of Mathematics,
University Walk, Bristol BS8 1TW, U.K.
\endaddress
\email v.ten\@bris.ac.uk
\endemail

\subjclass Primary 11D45, 11L07; Secondary 52B55
\endsubjclass

\keywords ``Hard'' Erd\"os problems, geometric complexity,
incidences, convex sequences, diophantine equations, exponential
sums, Falconer distance problem
\endkeywords

\thanks Research supported in part by the NSF grant
DMS02-45369, the Nuffield Foundation grant NAL/00485/A, and the
EPSRC grant GR/S13682/01.\newline
\endthanks
\date November 16, 2003
\enddate
\abstract We show that the equation
$$b_{i_1}+b_{i_2}+\dots+b_{i_d}=b_{i_{d+1}}+\dots+b_{i_{2d}} \tag*$$ has
$O\left(N^{2d-2+2^{-d+1}}\right)$ solutions for any strictly
convex sequence $\{b_i\}_{i=1}^N$ without any additional
arithmetic assumptions. The proof is based on weighted incidence
theory and an inductive procedure which allows us to effectively
deal with higher dimensional interactions. We also explain a
connection between this problem and the Falconer distance problem
in geometric measure theory.
\endabstract

\endtopmatter
\document
\loadeufm

\def \R{{\Bbb R}}
\def \Z{{\Bbb Z}}
\def \N {\frak N}

\head Section 1: Introduction and statement of results \endhead
\vskip.125in

Consider a sequence of real numbers ${\{b_i\}}_{i=1}^N$. It is a
classical problem in number theory to determine the number,
$\N_d=\N_d(N)$, of solutions of the equation
$$ b_{i_1}+b_{i_2}+\dots+b_{i_d}=b_{i_{d+1}}+\dots+b_{i_{2d}}.
\tag1.1$$

See, for example a book by Nathanson (\cite{Nath96}) or a survey
by Heath-Brown (\cite{HB02}), and the references contained therein
for a thorough description of various algebraic and combinatorial
aspects of this problems. The properties of $\N_d$ depend on
geometric and arithmetic properties of the sequence ${\{b_i\}}$.
For instance, if $b_i=i$, the number of solutions of $(1.1)$ is
$\approx N^{2d-1}$. \footnote{Here and throughout the paper the
notations $a \lesssim b$, or $a=O(b)$ means that there exists
$C>0$ such that $a \leq Cb$; the notation $a\gtrsim b$ means $b
\lesssim a$, and $a \approx b$ means that $a \lesssim b$ and $b
\lesssim a$. Similarly, $a \lessapprox b$ with a parameter $N$
means that for every $\epsilon>0$ there exists $C_{\epsilon}>0$
such that $a \leq C_{\epsilon}N^{\epsilon}b$.}

More interesting bounds are available if the sequence ${\{b_i\}}$
is a strictly convex in the sense that the points $\{(i,b_i)\}$
lie on a strictly convex curve in $\R^2$. For example, if
$b_i=i^2$, $\N_d \lesssim N^{2d-2}$ if $d \ge 4$, and the same
estimate with an appropriate power of $\log(N)$ if $d=2,3$. See,
for example the Landau's classical text (\cite{La69}) for a
thorough description of this and related issues. This example
shows that for a general strictly convex sequence, the best we can
hope for is an estimate of the form $\N_d \lessapprox N^{2d-2}$.
There are no examples known to the authors where
$\N_d>N^{2d-2+\epsilon}$ for $\epsilon>0$.

Taking $b_i=i^k$ puts us in the realm of the celebrated Waring
problem. It is conjectured that if $d=2$ and $k \ge 5$, the
equation $(1.1)$ only has trivial solutions. Similarly, in higher
dimensions it is conjectured that if $d$ is fixed and $k$ is
sufficiently large, then $(1.1)$ holds only if
$b_{i_l}=b_{i_{l+d}}$, up to permutation. See \cite{HB02} and the
references contained therein. These conjectures and the known
positive results in this context show that the bound $\N_d
\lessapprox N^{2d-2}$ can be improved under additional assumptions
on the arithmetic structure of the sequence $\{b_i\}$. Another
natural, non-integer example illustrating the power of arithmetic
considerations can be constructed as follows. Let $\{k_j\}$ be a
sequence of square free positive integers, such such that the
sequence given by $b_j=\sqrt{k_j}$ consists of numbers, such that
any subset of $d$ elements thereof is linearly independent over
${\Bbb Q}$. Then one can check that $\N_d \approx N^d$. In the
case $d=2$, one checks that the minimal polynomial of
$\sqrt{k_1}+\sqrt{k_2}$ is
$x^4-2x^2(k_1^2+k_2^2)-{(k_1^2-k_2^2)}^2$ and the assertion
follows by comparing the coefficients using the uniqueness of the
minimal polynomial. Indeed, it turns out that if
$b_{i_1}+b_{i_2}=b_{i_3}+b_{i_4}$, then
$(b_{i_1},b_{i_2})=(b_{i_3},b_{i_4})$, up to permutations. For
$d>2$, one uses the fact that the coefficients of the minimal
polynomial are symmetric polynomials of the roots. One then
observes that linear combination of the roots can be uniquely (up
to permutation) realized via these symmetric polynomials.

The main thrust of this paper is to obtain the best possible bound
on $\N_d$ under the assumption of strict convexity without any
additional arithmetic assumptions. It is reasonable to conjecture
that for every strictly convex sequence $\{b_i\}$, $\N_d
\lessapprox N^{2d-2}$. We prove that this estimate is
asymptotically true with an exponentially vanishing error as $d
\to \infty$. More precisely, we show (see Theorem 1 below) that
$$\N_d \lesssim N^{2d-2+2^{-d+1}}. \tag1.2$$

Konyagin (\cite{Ko03}) proved \footnote{This estimate is also
implicit in the results of the paper of Elekes et al.,
\cite{ENR99}.} that $(1.2)$ holds in the two-dimensional case. His
approach is based on the Szemer\'edi-Trotter incidence theorem.
When $d>2$, one is naturally led to consider issues associated
with higher dimensional incidence theory where significant
complication immediately arise. In this paper we address these
complication by a introducing an appropriate weighted version of
the Szemer\'edi-Trotter incidence theorem and an inductive
approach which we believe will find other interesting applications
in combinatorial geometry and additive number theory.

\subhead Statement of results \endsubhead Fix a convex sequence
${\{b_i\}}_{i=1}^N$, $N$ large, let ${\Cal N}\equiv \{1,2, \dots,
N\}$), and let $f:{\Bbb R} \to {\Bbb R}$ be a fixed strictly
convex function such that $f(i)=b_i$. Let ${\Cal
B}=\{b_1,\ldots,b_N\}$. Uniformity of the ensuing estimates is
understood in the sense that that none of the constants, hidden in
the estimates, and different for different $d$, depend on the
specific sequence ${\{b_i\}}_{i\in{\Cal N}}$ or $f$.

We need the following notation. The bounds for the quantities
$\N_d$ can be obtained by studying the set $${\Cal
C}_d\;\equiv\;\underbrace{{\Cal B}+\ldots+{\Cal B}}_{d \text{
times}}\;=\;\{c|\,c=b_{i_1}+\ldots+b_{i_d},\,\forall(i_1,\ldots,i_d)
\in{\Cal N}^d\}.\tag1.3$$

More precisely, for some $c\in {\Cal C}_d$ we shall refer to the
quantity
$$ \nu(c)\;= \;{|\{i_1,\ldots,i_d)\in {\Cal
N}^d:\,b_{i_1}+\ldots+b_{i_d}=c\}|\over d!}\tag1.4$$ as the {\it
weight} of $c$. $ $ \footnote{In $(1.4)$ above and throughout the
paper the notation $|\cdot|$ denotes the cardinality of a (finite)
set $\cdot$.}

We have
$$ \sum_{c\in {\Cal C}_d}\nu(c)={1\over d!}N^d, \tag1.5$$ the
right-hand side being the {\it net weight}, and
$$ \N_d=\sum_{c\in {\Cal C}_d}\nu^2(c). \tag1.6$$

Suppose the set ${\Cal C}_d=\{c_1,c_2,\ldots,c_t,\ldots\}$ has
been ordered by non-increasing weight. We shall see that in order
to estimate $\N_d$, it is sufficient to estimate the minimum
cardinality $|{\Cal C}_d|$ along with a majorant for the {\it
weight distribution} function, i.e. a decreasing function $\frak
n(t)$ such that $\frak n(t)\geq \nu(c_t)$. The inverse, also
decreasing function $\frak n^{-1}$ would provide the
bound\footnote{Note that $\frak n^{-1}$ is simply the distribution
function for $\frak n$ in the measure-theoretical sense.}
$$\frak n^{-1}(s)\;\geq\;| \ {\Cal C}_{d,s}\equiv\{c\in {\Cal
C}_d:\,\nu(c)\geq s \}|.\tag1.7$$

Our main result is the following.

\proclaim{Theorem 1} For $d\geq2$, let
$\alpha=2(1-2^{-d}),\,\beta=d-{4\over3}(1-2^{-d})$. Then
$$\hskip-.25in|{\Cal C}_d|\;\gtrsim \;N^{\alpha},\tag1.8$$
$$\hskip.05in\frak n(t) \;\lesssim \;N^{\beta} t^{-1/3},\tag1.9$$
$$\N_d\;\lesssim \;N^{2d-\alpha}.\tag1.10$$
\endproclaim

\remark{Remark} The main estimates of the theorem are $(1.9)$ and
$(1.10)$,  the latter one being the estimate on the number of
solutions of the diophantine equation $(1.1)$. The estimate
$(1.8)$ on cardinality of the sumset ${\Cal C}_d$ has been
included in the statement for the sake of completeness. This
estimate is implicit in \cite{ENR99}, Ch. 4 and is based on the
repeated application of the classical Szemer\'edi-Trotter theorem.
It is insufficient, however, to obtain estimates $(1.9)$ and
$(1.10)$, dealing the weight distribution within the set ${\Cal
C}_d$. This necessitated the development of certain weighted
incidence estimates included in the sequel.

Moreover, the proof of Theorem 1 applies to a more general class
of problems, such as counting the number of integer solutions of
the equation of the form $F(i_1, \dots, i_d)=F(i_{d+1}, \dots,
i_{2d})$, where $F$ is a strictly convex function of $d$ variables
which satisfies additional, fairly mild, assumptions on its
coordinate lower dimensional sections. \endremark

\vskip.1in It is interesting to contrast Theorem 1 with the
following well-known result from additive number theory, due to
Freiman (\cite{Frei73}).

\proclaim{Freiman's theorem} Let ${\Cal B} \subset {\Bbb Z}$ have
cardinality $N$, and suppose that $|{\Cal B}+{\Cal B}| \leq CN$.
Then ${\Cal B}$ is contained in a proper $s$-dimensional
progression\footnote{The set $P=\left\{x_0+\sum_{j=1}^s \lambda_j
x_j: 0 \leq \lambda_j < l_j \right\}$ with $x_0, \dots, x_s \in
{\Bbb Z}$ and $\lambda_1,\dots,\lambda_s;\,l_1, \dots, l_s \in
{\Bbb Z}^{+}$, is said to be an $s$-dimensional arithmetic
progression of length $l=\prod_{j=1}^s l_j$. $P$ is proper if
$l=|P|.$ } $P$ of length at most $KN$, where $s$ and $K$ depend
only on $C$.
\endproclaim

The estimate $(1.8)$ gives us $|{\Cal B}+{\Cal B}|\gtrsim
N^{\frac{3}{2}}$, and it is conjectured that the right bound is
$\gtrapprox N^2$. From the point of view of Freiman's theorem, a
power estimate is reasonable: if ${\{b_j\}}_{j=1}^N$ happens to be
a sequence of integers, the strict convexity assumption guarantees
that ${\Cal B}$ is not contained in an $s$-dimensional arithmetic
progression. However, a tighter connection with Freiman's theorem,
explaining the above exponent would be quite valuable. See
\cite{Gr02} for a description of Freiman's theorem and related
results on the structure of sumsets. See also \cite{KT99} and
\cite{KT01} for the description of related ideas in the context of
the Kakeya problem and the Falconer conjecture. Also see
\cite{Bo01} for the description of related issues in the context
of $\Lambda_p$ sets.

In Section 5 below we give a proof of a weaker though more robust
version of $(1.8)$ using Fourier analysis and results related to
the Falconer distance problem. This approach allows one to obtain
estimates on the size of ${\Cal C}_d$ with an additional
restriction that elements of this set be separated on a scale
depending on $N$. See \cite{HI2003}, \cite{IL2003}, and
\cite{IL2004}, where similar connections are explored. See also
\cite{Mag02} for a connection between diophantine equations and
ergodic theory.

In the case when the sequence $\{b_i\}_{i\in{\Cal N}}$ is
integer-valued, the estimate $(1.10)$ implies an estimate for the
$L_p$-norm of trigonometric polynomials with frequencies in
$\{b_i\}_{i\in{\Cal N}}$, i.e the Dirichlet kernel associated with
the sequence $\{b_i\}$.

\proclaim{Corollary 2} If $\{b_i\}_{i\in{\Cal N}}\subset\Z$, let
$$ f_N(\theta)=\sum_{j=1}^N e^{2 \pi i b_j \theta}.$$
Then
$$ \|f_N\|_{2d}\;\equiv\;{\left(\int_0^{2 \pi} {|f_N(\theta)|}^{2d}
d\theta \right)}^{\frac{1}{2d}} \;=\;O\left( N^{1-{1-2^{-d}\over
d}}\right).\tag1.11$$
\endproclaim
\remark{Remark} By expanding the square we see that $(1.10)$ and
$(1.11)$ are essentially identities when $d=1$. When $d>1$ observe
that $(1.11)$ is much stronger than the estimate that can be
obtained by interpolating the case $d=1$ and $d=\infty$ using
Holder's inequality. \endremark

\vskip.125in

\head{Section 2: Incidence theorems}\endhead

\vskip.125in

As we mention in the introduction, the main tool used in
\cite{ENR99} and \cite{Ko03} is the theorem of Szemer\'edi and
Trotter (\cite{ST83}) bounding the number of incidences between a
collection of points and straight lines in the Euclidean plane.
The theorem was extended to the case of points and hyper-planes or
spheres (with some natural restrictions on the arrangements) by
Clarkson et al. (\cite{CEGSW90}), see also the references therein.
It provides a powerful tool for solving problems in geometric
combinatorics. See also the books by Pach and Agarwal
(\cite{PA95}) and Matou$\check{\text s}$ek (\cite{Ma02}) for an
exhaustive description of this subject and related issues. It was
observed by Sz\'ekely (\cite{Sz97}) that the geometric graph
theory can deliver a short formal proof of the following statement
of the Szemer\'edi-Trotter incidence theorem in dimension two,
with the set of lines generalized to a class of curves satisfying
generic intersection hypotheses. \footnote{There is nothing to
prevent one from generalizing the ambient space $\R^2$ to a
general two-manifold of finite genus.} From this point on, we
shall use the terms ``lines'' and ``curves'' interchangeably.

\proclaim{Theorem 3 [Szemer\'edi-Trotter, Sz\'ekely]} Let $({\Cal
L},{\Cal P})$ be an arrangement\footnote{By the arrangement we
further mean an embedding, or drawing of the curves and points in
the plane.} of $m$ curves and  $n$ points in $\R^2$. Suppose that
no more than $\mu$ curves pass through any pair of points of
${\Cal P}$ and that any two curves of ${\Cal L}$ intersect at no
more than $\nu$ points of ${\Cal P}$. Then the total number of
incidences $$ I=|\{(l,p) \in {\Cal L} \times {\Cal P}: p \in l\}|
\;\lesssim\; {(\mu \nu)}^{\frac{1}{3}} {(mn)}^{\frac{2}{3}}+m+\mu
n. \tag2.1$$
\endproclaim

In the case of points and straight lines, $\mu=\nu=1$. Let us
further refer to this as the {\it simple intersection} case, where
$$ I\;\lesssim\; {(mn)}^{\frac{2}{3}}+m+n. \tag2.1a$$

The quantities $\mu$ and $\nu$, if bounded independently of  $m$
and $n$, may be viewed as constants which get absorbed into the
$\lesssim$ signs. Consequently, the assumption that any two curves
intersect at a {\it finite} (i.e. independent of $m,n$) number of
points, and that through any two points there pass no more than a
{\it finite} number of curves mean essentially that one has the
estimate given by $(2.1a)$ rather than the one given by $(2.1)$.

In the simple intersection case, the number of incidences $I$ for
the arrangement $({\Cal L}, {\Cal P})$ can be expressed in terms
of the counting function $\delta_{lp}$ for the arrangement. More
precisely,
$$I=\sum_{l\in{\Cal L},p\in{\Cal P}}\delta_{lp},$$ where
$\delta_{lp}=1$ if $p \in l$, and $0$ otherwise.

For our applications it will be useful to give an essentially
equivalent formulation of Theorem 3. We shall refer to $\mu$ and
$\nu$ as maximum weights. The numbers $m$ and $n$ shall will
referred to as net weights.

For the remainder of the paper, let us always consider the simple
intersection case. However, given some $\mu,\nu,$ we assign to
each line $l\in{\Cal L}$ and each point $p\in{\Cal P}$ the weights
$\mu(l)\in \{1,\ldots,\mu\}$ and $\nu(p)\in\{1,\ldots,\nu\}$,
respectively (although the individual weights certainly don't have
to be integer), so that
$$\sum_{l\in{\Cal L}}\mu(p)=m,\;\;\sum_{p\in{\Cal P}} \nu(p)=n.$$

Let us call such a weight assignment a {\it weight distribution
with the maximum weights} $(\mu,\nu)$ {\it and the net weights}
$(m,n)$. A single pair $(l,p)\in{\Cal L}\times{\Cal P}$ will have
the weight $w_{lp}=\mu(l)\nu(p)\delta_{lp}$. Let the number of
weighted incidences be defined by

$$I\equiv\sum_{l\in{\Cal L},p\in{\Cal P}}w_{lp}. \tag2.2$$

A weighted version of Theorem 3 is formulated as follows.

\proclaim{Theorem 3a} For all simple intersection arrangements
$({\Cal L},{\Cal P})$ with net weight $(m,n)$, for all weight
distributions with maximum weights $(\mu,\nu)$, one has
$$I \lesssim{(\mu
\nu)}^{\frac{1}{3}} {(mn)}^{\frac{2}{3}}+ \nu m+\mu n.\tag2.3$$
\endproclaim

Theorem 3a is shown to straightforwardly follow from Theorem 3 in
the Appendix, after a simple weight rearrangement argument. Note
that for the right hand side of $(2.3)$ one has
$$ {(\mu \nu)}^{\frac{1}{3}} {(mn)}^{\frac{2}{3}}+\nu m+\mu
n=\mu\nu\left[
{\left({m\over\mu}{n\over\nu}\right)}^{\frac{2}{3}}+{m\over\mu}+
{n\over\nu}\right],\tag2.4$$ which suggests that the maximum
number of weighted incidences is achieved when there are
${m\over\mu}$ lines and ${n\over\nu}$ points with uniformly
distributed weights. \vskip.125in

Note that unless the weights are distributed uniformly, neither
$|{\Cal L}|$, nor $|{\Cal P}|$ enter the estimate $(2.3)$.
Suppose, for example that it is known that $|{\Cal L}|\gg
m\mu^{-1}$, that is the majority of the lines have weights,
smaller than $\mu$. Can one use a divide-and-conquer approach to
take advantage of the average weight $\bar\mu= {m\over |{\Cal
L}|}$ in the formula (2.3)? As far as the equation $(1.1)$ is
concerned, the answer is yes. Ii is stated by Lemma 6, which is
central for the proof of Theorem 1. Note that the maximum weight
for the elements of ${\Cal C}_d$ can be bounded trivially by
$N^{d-1}$, or, less trivially, by $N^{d{d-1\over d+1}}$ using the
following theorem of Andrews (\cite{An63}) (see also \cite{BL98}).

\proclaim{Theorem 4 [Andrews]} The number of vertices of a convex
lattice polytope\footnote{A lattice polytope is a polytope with
vertices in the integer lattice $\Z^d$.}  in $\R^d$ of volume $V$
is $O\left( V^{{d-1\over d+1}}\right)$.
\endproclaim

\remark{Remark} The proof of Theorem 1 below, driven by the
weighted incidence technology, also yields an upper bound for the
number of solutions of the following generalization of the problem
considered in this paper:
$$ b_{i_1}+b_{i_2}+\ldots+b_{i_d} =a_{j_1}+a_{j_2}+\ldots +
a_{j_{D}} $$ of equation $(1.1)$, where $\{a_i\}_{i\in{\Cal N}}$
is another convex sequence and $2\leq d\leq D$.
\endremark

\vskip.125in

\head Section 3: Proof of Theorem 1 \endhead

\vskip.125in

The proof is by induction on $d$, starting from $d=2$. Let
$$\gamma=\{(x, f(x)): x \in [1,N]\},
\;\;\text{and}\;\;\gamma_{\Cal N}=\{(i, f(i)): i \in {\Cal
N}\}\tag3.1$$ \specialhead

The case $d=2$:\endspecialhead

\proclaim{Lemma 5} We have
$$|{\Cal C}_2|\;\gtrsim\; N^{3/2},\tag3.2$$
and
$$ |{\Cal C}_s|=|\{c\in{\Cal C}_2:\,\nu(c)\geq s\}|\;\lesssim\;
N^3s^{-3}.\tag3.3$$
\endproclaim

\noindent{\bf Proof.} Define
$${\Cal N}_2\;\equiv {\Cal N}+{\Cal N}.$$
Consider the set of points ${\Cal P}={\Cal N}_2\times{\Cal
B}\,+\,\gamma_{\Cal N}$ and the set of curves ${\Cal
L}=\gamma+{\Cal N}_2\times{\Cal B}.$ Convexity implies that the
arrangement $({\Cal L},{\Cal P})$ satisfies the simple
intersection condition.

Since $|{\Cal P}|\lesssim N^2,$ the number of incidences $I$ for
this arrangement can be estimated by the non-linear term in
formula $(2.1a)$, i.e $$ I\;\lesssim\; N^{4/3} (|{\Cal P}|)^{2/3}.
\tag 3.4$$

On the other hand, each curve of ${\Cal L}$ contains at least $N$
points of ${\Cal P}$ (that is why ${\Cal P}$ has been taken as
${\Cal N}_2\times{\Cal B}$ rather than simply ${\Cal N}\times{\Cal
B}$). It follows that $I\gtrsim N^3$, and
$$ N|{\Cal C}_2| \;\approx \;|{\Cal P}|\; \gtrsim\;N^{5/2},$$ which implies
$(3.2)$.

Let ${\Cal P}_s=\{p\in{\Cal P}:\,\nu(p)\geq s\}$, where $\nu(p)$
is the number of curves of the arrangement ${\Cal L}$ intersecting
at the point $p$. Applying the estimate $(3.4)$ for the number of
incidences for the arrangement $({\Cal L},{\Cal P}_s)$, with
$|{\Cal P}_s|$ substituting ${\Cal P}$ and comparing it with the
lower bound $s|{\Cal P}_s|$, we get
$$|{\Cal C}_s|\;\approx\;N^{-1}|{\Cal P}_s|\;\lesssim\;N^3 s^{-3},$$
which implies $(3.3)$.

\vskip.1in

In view of $(3.2)$ let $\bar\nu=\sqrt{N}$ be the (approximate)
upper bound for the average weight per element of ${\Cal C}_2$. By
$(3.3)$ the weight distribution function in the (ordered) set
${\Cal C}_2$ satisfies
$$ \nu(c_t)\;\lesssim\;\frak n(t) \;=\; Nt^{-1/3}.$$

It follows that for the set ${\Cal C}_{2,\bar\nu},$ containing
those $O(N^{3/2})$ elements of ${\Cal C}_2,$ whose weights may
exceed $\bar\nu$, one has
$$\sum_{c\in {\Cal C}_{2,\bar\nu}} \nu(c)^2\lesssim N^2
\int_1^{N^{3/2}}t^{-2/3}dt\;\approx\;N^{5/2}.\tag3.5$$ On the
other hand, for the complement ${\Cal C}_{2,\bar\nu}^c$ of ${\Cal
C}_{2,\bar\nu}$ in ${\Cal C}_2$, where the weight does not exceed
$\bar\nu$, one has $$\sum_{c\in {\Cal C}^c_{2,\bar\nu}}
\nu(c)^2\lesssim \bar\nu \sum_{c\in {\Cal C}_2}\nu(c)\;\approx
\;N^{5/2},\tag3.6$$ as the total weight of ${\Cal C}_2$ is
approximately $N^2$. This proves the formulas $(1.8--1.10)$ for
$d=2$.

\remark{Remark} Formulas $(3.5)$ and  $(3.6)$ are easily
understood in the sense of the defining formulas $(1.4-1.7)$,
dealing with the weight distribution function $\nu(c)$ in the set
${\Cal C}_2$, with the known $L_1$ norm of $\nu$, the net weight.
The quantity ${\N}_2$ in question, which is the square of the
$L_2$ norm of $\nu(c)$ can be bounded as follows. One naturally
partitions the domain ${\Cal C}_2$ in two subsets. In the first
subset, where $\nu(c)$ is likely to exceed an upper bound for its
average $\bar\nu$ (obtained as the net weight divided by the lower
bound for $|{\Cal C}_2|$) one uses the (strictly decreasing,
concave) majorant $\frak n(t)$ for $\nu(c_t)$ and gets $(3.5)$.
The integral of $\nu^2(c)$ over the second subset, where
$\nu(c)\lesssim \bar\nu$ is bounded by the product of the $L_1$
norm of the function $\nu(c)$ ($\approx N^2$) and the $L_\infty$
norm $\bar\nu=\sqrt{N}$ over the subset: this is $(3.6)$. The same
tactics is used in the following main part of the proof. The
tricky part there is getting the tight enough majorant $\frak
n(t).$
\endremark

\specialhead The case $d\Rightarrow d+1$:\endspecialhead In order
to characterize the weight distribution function $\nu(c)$, for
$c\in{\Cal C}_{d+1}$, consider the equation $$ f(i_1) \,+
\,[\,f(i_2)+\ldots+f(i_{d+1})]\, =\, c. \tag3.7$$ Let $a\in {\Cal
C}_d$. Extend $(3.7)$ to the system of equations
$$ \left\{\matrix
f(i_1) + u = c, \\ \;\;\;\;\;i_1\;+ j= k,
\endmatrix\right.
\qquad\forall\,(i_1,j,k,u,c)\in{\Cal N}\times{\Cal N}_2\times{\Cal
N}_2\times {\Cal C}_d\times{\Cal C}_{d+1}.\tag3.8$$

Note that ${\Cal C}_{d+1}$ is considered as a set, rather than
multi-set. The elements of the set ${\Cal
C}_d=\{u_1,u_2,\ldots,u_t,\ldots\}$  are endowed with
non-increasing weights, with some weight distribution function
$\mu(u)$. In particular,  the $L_1$ norm of $\mu(u)$, $\|\mu\|_1$,
over ${\Cal C}_d$ is $O(N^d)$, the $L_\infty$ norm is
$O(N^{d{d-1\over d+1}}),$ by the aforementioned Andrews theorem,
and (by the induction assumption) that there is a majorant
$$\mu(u_t)\;\lesssim\;\frak
m(t)\;=\;N^{\beta_d}t^{-1/3},\tag3.9$$ where $\beta_d=d
-{4\over3}(1-2^{-d})$. There is also the estimate $(1.8)$ for the
cardinality of ${\Cal C}_d$, enabling one to introduce the upper
bound for the average weight $\bar\mu$ in ${\Cal C}_d$ as follows:
$$ {\|\mu\|_1\over|{\Cal C}_d|}\;\lesssim\;\bar\mu\;=\;
N^{d-\alpha_d},\tag3.10 $$ with $\alpha_d=2-2^{-d+1}$.

The number of solutions of $(3.7)$ is not smaller than the number
of solutions of $(3.8)$, divided by $N$. The number of solutions
of $(3.8)$ can be estimated in terms of the number of weighted
incidences $I$ between the weighted set ${\Cal L}$ of the curves,
given by the translations $\gamma_{ju}$ of the curve $\gamma$
defined by $(3.1)$, by the elements of ${\Cal N}_2\times{\Cal
C}_d$ and the set ${\Cal P}={\Cal N}_2\times{\Cal C}_{d+1}$. Thus
the problem essentially boils down to the same scheme as it was in
the case $d=2$, except that {\it weighted} incidences should be
counted in order to verify estimates $(1.9)$ and $(1.10)$.
Verification of $(1.8)$ is easier: it requires only the available
(through the induction assumption) lower bound $|{\Cal
C}_d|\gtrsim N^{\alpha_d}$ and the use of $(2.1a)$ and was done in
\cite{ENR99} (and in the case $d=2$, see $(3.4)$ and the formula
that follows it). The corresponding estimate $(3.12)$ can be also
obtained via the ensuing Lemma 6, which we have chosen to do in
order to show that the lemma by itself is tight.

In our consideration, each translated curve $\gamma_{ju}$ would
inherit the weight $\mu(u)$ of the corresponding element $u\in
{\Cal C}_d$. The following lemma is central for the rest of the
proof.

\proclaim{Lemma 6} Under the assumptions $(3.9)$ and $(3.10)$ on
the weight distribution function $\mu(u)$ in the set ${\Cal C}_d$,
the number of incidences for the above defined arrangement $({\Cal
L},{\Cal P})$, describing the solutions of the system $(3.8)$ is
bounded as follows:
$$ I\;\lesssim\; \bar\mu^{1/3} N^{2(d+1)/3}(N|{\Cal
C}_{d+1}|)^{2/3}.\tag3.11$$
\endproclaim

Lemma 6 shows that in order to count the weighted incidences in
the arrangement $({\Cal L},{\Cal P})$, instead of the maximum
weight upper bound $O(N^{d{d-1\over d+1}})$ in the set ${\Cal L}$
(transferred from ${\Cal C}_d$), given by the Andrews theorem, one
can use the formula $(2.3)$ with the (smaller) average weight
majorant $\bar\mu$ for the quantity $\mu$, as well as (naturally)
the net weight $m=N^{d+1}$ and the maximum point weight $\nu=1$ in
${\Cal P}$. The proof of Lemma 6 is given in the next section. We
shall now use it to complete the proof of Theorem 1.

\vskip.1in

Assuming Lemma 6, we compare its estimate $(3.11)$ with the fact
that on each curve of ${\Cal L}$ there lies at least $N$ points of
${\Cal P}$, thus $I\geq N^{d+2}$, because $N^{d+1}$ is
approximately the net weight of ${\Cal L}$. Comparing the powers
of $N$, we get
$$|{\Cal C}_{d+1}|\;\gtrsim\; N^{2-2^{-d}}\,=\,N^{\alpha_{d+1}}.\tag3.12$$

This leads us to define the upper bound for the average weight in
${\Cal C}_{d+1}$
$$ \bar\nu\;=\;N^{d+1-\alpha_{d+1}}.\tag3.13 $$

Let ${\Cal P}_s=\{p\in{\Cal P}:\,\nu(p)\geq s\}$, where $\nu(p)$
is now defined as the total weight of all the curves of the
arrangement ${\Cal L}$ intersecting at the point $p$. Clearly
${\Cal P}_s={\Cal N}_2\times {\Cal C}_{d+1,s}$, where ${\Cal
C}_{d+1,s}$ is the subset of ${\Cal C}_{d+1}$, consisting of all
those elements whose weight is not smaller than $s$. In order to
estimate $|{\Cal C}_{d+1,s}|$, formula $(2.1a)$ cannot be used, as
one has to take into account the individual weight of each curve
$\gamma_{ju}\in {\Cal L}$, passing through the given point $p$.
Instead, weighted incidences have to be dealt with, and Lemma 6
enables one use the average weight $\bar\mu$ in the estimate,
rather than the maximum weight $\mu\gg\bar\mu$.

In view of this, we proceed by comparing the lower bound $s N
|{\Cal C}_{d+1,s}|$, for the number of weighted incidences for the
arrangement $({\Cal L},{\Cal P}_s)$ with $(3.11)$, in which
$|{\Cal C}_{d+1,s}|$ substitutes ${\Cal C}_{d+1}$. This yields
$$|{\Cal C}_{d+1,s}|\;\lesssim\;N^{1-\alpha_d}\left({N^d\over
s}\right)^3.\tag3.14$$ If $s=\bar\nu$, defined by $(3.10)$, it
follows that

$$|{\Cal C}_{d+1,\bar\nu}|\;\lesssim \;N^{\alpha_{d+1}},\tag3.15$$
which is the same as the right-hand side in $(3.12)$, and complies
with $(1.8)$. Inversion of $(3.14)$ yields the majorant for the
weight distribution function $\nu(c)$ for $c\in{\Cal C}_{d+1}$:
$$ \nu(c_t)\;\lesssim\;\frak n(t)\;=\; N^{d-{1-2^{-d+1}\over3}}
t^{-1/3}\;=\;N^{\beta_{d+1}}t^{-1/3}, \tag3.16$$ as is claimed by
$(1.9)$.

The final step of the proof follows the remark at the end of the
$d=2$ section. Namely one partitions
$${\Cal C}_{d+1}\;=\;{\Cal C}_{d+1,\bar\nu}\;\cup\;{\Cal
C}_{d+1,\bar\nu}^c,$$ the first piece containing ``heavy''
elements, and estimates
$$ \sum_{c\in{\Cal C}_{d+1,\bar\nu}^c}
\nu^2(c)\;\lesssim\;N^{d+1}\bar\nu\;=\;N^{2(d+1)-\alpha_{d+1}},\tag3.17$$
as well as
$$\sum_{c\in{\Cal C}_{d+1,\bar\nu}}
\nu^2(c)\;\lesssim\;N^{2\beta_{d+1}}\int_1^{N^{\alpha_{d+1}}}t^{-2/3}dt\;
\approx\;N^{2(d+1)-\alpha_{d+1}}.\tag3.18 $$

The estimates $(3.17)$ and  $(3.18)$ are consistent with $(1.10)$.
Thus the proof of Theorem 1 is complete.

\vskip.125in

\head Section 4: Proof of Lemma 6 \endhead

\vskip.125in

The objective is to partition the set
$${\Cal C}_d\,=\,\bigcup_{i=0}^M{\Cal C}_{d,\mu_i} \tag4.1 $$ into
$M$ (a fairly large number of) pieces, trying to make each one of
them as large as possible, yet having control over the number of
weighted incidences it can possibly be responsible for. For
simplicity let
$$ {\Cal C}\;\equiv\;{\Cal C}_d,\;\,{\Cal C}_i\;\equiv\;{\Cal
C}_{d,\mu_i}. $$

The partition is required to have the following property:
$$\mu(c)\,\lesssim \,\mu_i,\;\,\forall c\in{\Cal
C}_{i},\tag4.2$$ where the strictly decreasing sequence $\mu_i$
will start out from
$$ \mu_0\,=\,N^{d{d-1\over d+1}}, $$ (the maximum weight granted
by the Andrews theorem \footnote{In fact, one can see from the
proof that the use of the Andrews theorem is unnecessary: one can
simply start out with $\mu_0=N^{d}$, which is the net weight of
${\Cal C}$.}) and go down geometrically to the average weight
$\bar\mu$ in ${\Cal C}$, specified in $(3.10)$. The number $M$ is
chosen in such a way that $\mu_M$ gets sufficiently close to
$\bar\mu$, so that the effect of the difference between them can
be swallowed by the $\lesssim$ symbol. The sequence $\{{\Cal
C}_{i}\}$ will be constructed, using the weight distribution
majorant $(3.9)$.

By the general estimate $(2.3)$ of Theorem 3a in order to prove
the lemma, it suffices to show that
$$ \left(\tilde I\;\equiv\;\sum_{i=0}^{M}  \mu_i^{1\over 3}
m_i^{2/3}\right)\;\lesssim\; \left(\bar I\;\equiv\;\bar \mu^{1/3}
m^{2/3}\right),\tag4.3 $$ where $m=N^d$ is the net weight of
${\Cal C}$, and $m_i$ is the net weight of each ${\Cal C}_i,$ for
$i=0,\ldots,M$. Indeed, it easy to see that the linear terms
coming from the bound $(2.3)$ are irrelevant: the first linear
term is $N^{d+1}$, being the total weight of the set of lines
${\Cal L}$, defined by the system of equations $(3.8)$; the second
linear term, relative to the set ${\Cal C}_i$ will be equal to
$\mu_i N^{d+1}$, as $N^{d+1}$, is also the net weight of the set
of points ${\Cal P}$, defined by the system of equations $(3.8)$.
Both terms will be dominated by the quantity $\tilde I$ defined by
$(4.3)$ by construction. This is shown explicitly in the end of
the proof.

The weights $m_i$ are to be estimated via $\mu_i$, using the
inverse formula for the majorant $(3.9)$, i.e
$$|\{c\in {\Cal C}:\,\mu(c)\geq s\}|\;\lesssim\;  \frak m^{-1}(s)\;=\;
N^{3\beta_{d}} s^{-3},\;\;\beta_d=d-{4\over3}(1-2^{-d}).\tag4.4$$
Note that the majorant $(3.9)$ is good for nothing as far as the
elements $c$ of ${\Cal C}$, such that $\mu(c)\lesssim\bar\nu$ are
concerned. Indeed, a calculation yields
$$ \int_{\bar\mu}^\infty \frak m^{-1}(s)\,ds\,\approx\,m,$$
where $m\approx N^{d}$ is the net weight of ${\Cal C}$.

Also for the terms in the sum in the right-hand side of $(4.3)$
denote
$$ \tilde I_i\;\equiv\;\mu_i^{1\over 3} m_i^{2/3}. $$

The sets ${\Cal C}_i$ and the number $M$ are to be chosen such
that
$$ \tilde I_i\;\lesssim\; N^{-\varepsilon}\bar
I,\;\;M\,\approx\,N^{\varepsilon},\tag4.5 $$ for some fixed small
positive number $\varepsilon$.

Let us describe the first step of the construction.  Let a number
$\delta_0$ be defined via $\mu_0=N^{\delta_0}\bar\mu$. Define the
weight $m_0$ of the set ${\Cal C}_0$ implicitly via $(4.5)$, i.e.
$$ \mu_0^{1/3} m_0^{2/3}\,\approx\,N^{-\varepsilon}\bar\mu^{1/3}
m^{2/3}, $$ which yields
$$ m_0=N^{-{1\over2}(3\varepsilon+\delta_0)}\,m.\tag4.6$$ Then the
weight of any element $c$ in the complement ${\Cal C}_0^c$ of
${\Cal C}_0$ in ${\Cal C}$ should be bounded from above by some
quantity $\mu_1$, which can be defined implicitly from $$
\int_{\mu_1}^\infty \frak m^{-1}(s)\,ds\,=\,m_0. $$ This yields
$$ \mu_1=\bar\mu N^{\delta_1},\;\delta_1={1\over4}(3\varepsilon+\delta_0).
\tag4.7$$

One can see that for $\varepsilon$ small enough, say
$\varepsilon={1\over9}\delta_0$, one has
$$\delta_1\leq {1\over3}\delta_0.$$

The procedure is now repeated for the set ${\Cal C}_0^c$, where
the maximum weight is bounded in terms of $\mu_1$, rather than
$\mu_0$, which will result in some set ${\Cal C}_1$ having been
pulled out of it, such that the maximum weight in the complement
of ${\Cal C}_1$ in ${\Cal C}_0^c$ is bounded in terms of some
$\mu_2\ll\mu_1$, and so on. After having done it $M-1$ times, the
set ${\Cal C}$ will be partitioned, according to $(4.1)$, where
the last member of the partition ${\Cal C}_M$ is the complement of
the union $\bigcup_{i=0}^{M-1}{\Cal C}_i$ in ${\Cal C}$. For
$i=1,\ldots,M$ the maximum individual element weight in ${\Cal
C}_i$ is bounded similarly to $(4.7)$, namely
$$\mu_i=\bar\mu\,N^{\delta_i},\;\,\delta_i={1\over4}(3\varepsilon+\delta_{i-1}).
\tag4.8$$

By construction, each set ${\Cal C}_i$, for $i=0,\ldots,M-1$ would
create the number of incidences $I_i$ with ${\Cal P}$, bounded (by
$(3.11)$ and $(4.5)$)  by
$$ I_i\;\lesssim\;N^{2d+2^{-d}-\varepsilon+1}. $$

Note that in comparison with $(1.10)$ one has $d\rightarrow d+1$,
which accounts for an extra $N$ here, as the quantity ${\N_d}$
equals $N^{-1}$ times the number of incidences for the arrangement
$({\Cal L},{\Cal P})$, introduced apropos of the system of
equations $(3.8)$, rather than equation $(3.7)$.

For  $\varepsilon\leq1$, the right hand side of the last
expression will exceed the maximum for the linear term in the
estimate $(2.3)$, applied to the arrangement $({\Cal L},{\Cal
P}),$ as the latter can be bounded simply via
$$\mu_0 N^{d+1}\;\lesssim\;N^{2d^2\over d+1}.$$

Finally, one is left to estimate
$$ \mu_M\,\leq \, \bar\mu
N^{{1\over3^{M}}\delta_0}\;\lesssim\;\bar\mu,$$ by $(4.7)$,
$(4.8)$. Indeed, $M=N^{\varepsilon}$, and we see that the
remaining set ${\Cal C}_M$ will not be responsible for more
incidences than specified by the right-hand side of $(3.11)$. This
completes the proof of Lemma 6.

\vskip.125in

\head Section 5: The Falconer distance problem and convex
sequences \endhead

\vskip.125in

In this section we recover some of the estimates implied by
Theorem 1 using Fourier analysis and connect the arithmetic
problem studied in this paper to the Falconer distance conjecture
in geometric measure theory.

The Falconer distance problem (see e.g. \cite{Fa86} asks whether
the Lebesgue measure ${\Cal L}^1$ of the distance set
$\Delta(E)=\{|x-y|: x,y \in E\}$, $|x|=\sqrt{x_1^2+\dots+x_d^2}$,
$d \ge 2$, is positive provided that the Hausdorff dimension $dim$
of $E \subset {[0,1]}^d$ is sufficiently large. It is conjectured
that the conclusion should hold provided that $ dim(E)\,>\,
\frac{d}{2},$ and there exists and an arithmetic example, based on
the integer lattice and diophantine approximation shows that such
a result would be the best one possible. Namely, Falconer
(\cite{Fa86}) proved the first result in this direction, having
shown that ${\Cal L}^1(\Delta(E))>0$ in ${\Bbb R}^d$, provided
that $dim(E)>\frac{d+1}{2}$.

The best known result in two dimensions is due to Wolff
(\cite{Wo99}) who proved that ${\Cal L}^1(\Delta(E))>0$ provided
that $dim(E)>\frac{4}{3}$. See also previous improvements due to
Bourgain (\cite{Bo94}). In higher dimensions, the best known
estimate is due to Erdogan (\cite{Er03}) who proved that ${\Cal
L}^1(\Delta(E))>0$ provided that
$$dim(E)>\frac{d(d+2)}{2(d+1)}.\tag5.1$$
Moreover, Erdogan's proof makes it clear that the Euclidean
distance in the definition of the distance set may be replaced by
any distance $\rho$ such that the level set $\{x: \rho(x)=1\}$ is
smooth, convex and has curvature, bounded from below.

Recall the definition of the set ${\Cal C}_d$ from $(1.3)$:
$c\in{\Cal C}_d$ should satisfy the equation
$$ b_{i_1}+b_{i_2}+\dots+b_{i_d}=c.$$
Let $f$ be the convex function, underlying the sequence $\{
b_i\}$, i.e. such that $$f(i)=b_{i}.\tag5.2$$

Let $q_1=2$, and $q_{j+1}=q_j^j$. Let
$$ E_j=\left\{x \in {[0,1]}^d: \left|x_k-\frac{p_k}{q_j}\right| \leq
q_j^{-\frac{d}{s}},\, \;\text{for every}\  p=(p_1,\dots,p_d) \in
{\Bbb Z}^d \cap {[0,q_j]}^d \right\}, \tag5.3$$ for some
$s\in(0,d)$.

Let $E=\cap E_j$. It follows from Theorem 8.17 in \cite{Fa85} that
the Hausdorff dimension of $E$ is $s$.

Define $\rho_f(x)=f(x_1)+\dots+f(x_d)$ and
$\Delta_{\rho_f}(E)=\{\rho(x-y): x,y \in E\}$. Then $f$ can be
chosen such that $(5.2)$ is satisfied and
$$ |\Delta_{\rho_f}(E_j)| \;\lesssim \;q_j^{-\frac{d}{s}} \times
|\{\rho_f(z-w): z,w \in {\Bbb Z}^d \cap {[0,q_j]}^d\}|. \tag5.4$$

Suppose that there exists a subsequence of the $q_j$s going to
infinity such that
$$
|\{\rho_f(z-w): z,w \in {\Bbb Z}^d \cap {[0,q_j]}^d\}|
\,\lesssim\, q_j^{\beta}. \tag5.5$$

Plugging this estimate into $(5.4)$ we see that
$$ |\Delta_{\rho_f}(E_j)| \,\lesssim \,q_j^{-\frac{d}{s}} q_j^{\beta}.
\tag5.6$$

Suppose that the function $f$ can be chosen such that $(5.2)$
holds, the level set
$$\{x: \rho_f(x)=1\}\tag5.7$$ is smooth and has
curvature bounded from below. It follows from $(5.1)$ that
$\Delta_{\rho_f}(E)$ has positive Lebesgue measure if
$s>\frac{d(d+2)}{2(d+1)}$. However, this is in direct
contradiction with $(5.6)$ if $\beta<\frac{2(d+1)}{d+2}$. In the
language of Theorem 1, we have just proved that
$$ |{\Cal C}_d|\, \gtrapprox \,N^{2-\frac{2}{d+2}}, \tag5.8$$ which
recovers $(1.8)$ of Theorem 1 in the case $d=2$. Observe that we
actually proved a little more, namely that the number of
$N^{-\frac{d-1}{d+1}}$-{\it separated values} of ${\Cal C}_d$ is
bounded from below by $N^{2-\frac{2}{d+2}}$, not just the total
number of values.

Note that the Falconer conjecture, if true would imply that
$$ |{\Cal C}_d|\, \gtrapprox \,N^{2},\;\forall d\geq2. \tag5.9$$

For $d\geq2$, the estimate $(5.7)$ is weaker that $(1.8)$, because
the method in this section clearly does not take any advantage of
the inductive procedure made possible by the special structure of
the function $\rho_f$. However, it is clear that the scope of this
method should be considerably wider, in the sense that it enables
one to consider more general arithmetic problems than those given
by sums of one-dimensional functions. In particular, the analytic
method in question can be applied in high generality to the study
of equation of the form
$$ F(i_1, \dots, i_d)=F(i_{d+1}, \dots, i_{2d}), \tag5.10$$ where $F$
is a general strictly convex function of $d$ variables. We shall
undertake a systematic study of these equations and connections
with arithmetic methods in a subsequent paper.

\vskip.125in

\head Appendix: Proof of Theorem 3a \endhead

\vskip.125in

Without loss of generality, one can assume that all the weights
are integers, the net line weight $m$ is a multiple of the maximum
line weight $\mu$, and the net point weight $n$ is a multiple of
the maximum point weight $\nu$. Then the formula $(2.4)$ is
equivalent to the bound $(2.1a)$, for the number of incidences
between $m/\mu$ lines and $n/\nu$ points, provided that each
incidence would be counted $\mu\nu$ times. In other words, for the
uniform weight distribution there is nothing to prove.

Otherwise, consider some arrangement $({\Cal L},{\Cal P})$ and
suppose, that the weight distribution over, say ${\Cal P}$ is not
uniform. Then there exist $p_1,p_2\in{\Cal P}$, such that
$\nu(p_1)<\nu(p_2)<\nu.$  For $p\in{\Cal P}$ let
$$w_p\,=\,\sum_{l\in{\Cal L}}\mu(l)\delta_{lp},$$ be the total
weight of all the lines incident to $p$. If $w(p_1)>w(p_2)$, first
change the weight distribution by swapping the values $\nu(p_1)$
and $\nu(p_2)$ over the points $p_1$ and $p_2$. Then modify the
weight distribution by changing $\nu(p_1)\rightarrow \nu(p_1)-1$
and $\nu(p_2)\rightarrow \nu(p_2)+1$. If $\nu(p_1)$ has become
zero, remove $p_1$ from ${\Cal P}$. As the result, the weight
distribution has been modified, so that the number of incidences
has increased, and the net weight has stayed constant. Continue
this (greedy) procedure, until the weight distribution over ${\Cal
P}$ has become uniform; then do the same thing with the set ${\Cal
L}$. At each single step, the number of incidences will have
increased. However, as the result, one still ends up with the
upper bound $(2.1a)$, as only $m/\mu$ lines and $n/\nu$ points
remain.

\newpage

\head References \endhead

\vskip.125in

\ref \key An63 \by G. E. Andrews \paper A lower bound for the
volume of strictly convex bodies with many boundary lattice points
\jour  Trans. Amer. Math. Soc. \yr 1963 \vol 106  \pages 270--279
\endref

\ref \key BL98 \by I. B\`ar\`any and D. G. Larman \paper The
convex hull of the integer points in a large ball. \jour  Math.
Ann. \yr 1998 \vol  312 \pages 167--181 \endref

\ref \key Bo94 \by J. Bourgain \paper Hausdorff dimension and
distance sets \jour Israel J. Math. \vol 87 \yr 1994 \pages
193--201 \endref

\ref \key Bo01 \by J. Bourgain \paper $\Lambda\sb p$-sets in
analysis: results, problems and related aspects \jour Handbook of
the geometry of Banach spaces; North-Holland, Amsterdam \vol I
\pages 195--232 \yr 2001 \endref

\ref \key CEGSW90 \by K. Clarkson, H. Edelsbrunner, L. Guibas, M.
Sharir, and E. Welzl \paper Combinatorial complexity bounds for
arrangements of curves and surfaces \jour Discrete and
Computational Geometry \yr 1990 \vol 5 \pages 99--160 \endref

\ref \key ENR99 \by G. Elekes, M. Nathanson, and I. Ruzsa \paper
Convexity and sumsets \jour Journal of Number Theory \yr 1999 \vol
83 \pages 194--201 \endref

\ref \key Er03 \by M. B. Erdogan \paper A note on the Fourier
transform of fractal measures \jour (preprint) \yr 2003
\endref

\ref \key Fa85 \by K. J. Falconer \paper The geometry of fractal
sets \jour Cambridge University Press \yr 1985 \endref

\ref \key Fa86 \by K. J. Falconer \paper On the Hausdorff
dimensions of distance sets \jour Mathematika \vol 32 \pages
206--212 \yr 1986 \endref

\ref \key Frei73 \by G. Freiman \paper Foundations of a structural
theory of set addition \yr 1973 \jour Translations of Mathematical
Monographs; AMS Providence, R.I. \vol 37 \endref

\ref \key Gr02 \by B. Green \paper Structure theory of set
addition \yr 2002 \jour Edinburgh lecture notes (available at
http://www.dpmms.cam.ac.uk/$\thicksim$bjg23/) \endref

\ref \key HB02 \by D.R. Heath-Brown \paper Counting rational
points on algebraic varieties \yr 2002 \jour (preprint) \endref

\ref \key HI03 \by S. Hofmann and A. Iosevich \paper Circular
averages and Falconer/Erdos distance conjecture in the plane for
random metrics \jour (submitted) \yr 2003 \endref

\ref \key Io02 \by Iosevich, A. \paper Szemer\'edi-Trotter
incidence theorem, related results and some amusing consequences
\jour to appear \endref

\ref \key IL2003 \by A. Iosevich and I. Laba \paper $K$-distance
sets and Falconer conjecture \jour (submitted for publication) \yr
2002 \endref

\ref \key IL2004 \by A. Iosevich and I. Laba \paper Cones over
varieties, generalized distance sets and representations of
integers \jour (in preparation) \yr 2004 \endref

\ref \key KT99 \by N. Katz and T. Tao \paper Bounds on arithmetic
progression, and applications to the Kakeya conjecture \jour Math.
Res. Let. \vol 6 \yr 1999 \pages 625--630 \endref

\ref \key KT01 \by N. Katz and T. Tao \paper Some connections
between Falconer's distance set conjecture and sets of Furstenburg
type \jour New York J. Math. (electronic) \vol 7  \yr 2001 \pages
149--187 \endref

\ref \key Ko03 \by S. Konyagin \paper (personal communication) \yr
2003 \endref

\ref \key La69 \by E. Landau \paper Vorlesungen \"Uber
Zahlentheorie (German) \jour Chelsea Publishing Co., New York \yr
1969 \endref

\ref \key Mag02 \by A. Magyar \paper Diophantine equations and
ergodic theorems \jour Amer. J. Math.  \vol 124  \yr 2002 \pages
921--953 \endref

\ref \key Ma02 \by J. Matou$\check{\text s}$ek \book Lectures on
Discrete Geometry \yr2002\bookinfo  Springer-Verlag, New
York\endref

\ref \key Nath96 \by M. Nathanson \book Additive Number Theory-
the classical bases \yr 1996 \bookinfo Springer \endref

\ref \key PA95 \by J. Pach and P. Agarwal \book Combinatorial
Geometry \yr 1995 \bookinfo Wiley-Interscience Series in Discrete
Mathematics and Optimization \endref

\ref \key Sz97 \by L. Sz\'ekely \paper Crossing numbers and hard
Erdos problems in discrete geometry \jour Combinatorics,
Probability, and Computing \yr 1997 \vol 6 \pages 353--358 \endref

\ref \key ST83 \by E. Szemer\'edi and W. Trotter \paper Extremal
problems in discrete geometry \jour Combinatorica  \vol 3 \yr 1983
\pages 381--392 \endref

\ref \key Wo99 \by T. Wolff \paper Decay of circular means of
Fourier transforms of measures \jour Internat. Math. Res. Notices
\yr 1999 \vol 10 \pages 547--567 \endref

\enddocument